\documentclass[11pt,reqno]{amsart}

\usepackage[margin=1.15in]{geometry}
\usepackage{amsmath,amssymb,amsthm,mathtools}
\usepackage[expansion=false]{microtype}
\usepackage{needspace}
\usepackage[hidelinks]{hyperref}

\newcommand{\N}{\mathbb N}
\newcommand{\e}{\mathrm e}
\newcommand{\Bset}{\mathcal B}
\newcommand{\ind}{\mathbf 1}

\theoremstyle{plain}
\newtheorem{theorem}{Theorem}[section]
\newtheorem{proposition}[theorem]{Proposition}
\newtheorem{lemma}[theorem]{Lemma}

\theoremstyle{definition}

\theoremstyle{remark}
\newtheorem{remark}[theorem]{Remark}

\title{A Counterexample to Wooley's
\texorpdfstring{$L^1$}{L1} Question for Beatty Sequences}

\author{Dijia Chen}
\address{Department of Mathematics,
University of Wisconsin--Madison,
Madison, Wisconsin, USA}
\email{dchen426@wisc.edu}



\begin{document}
\raggedbottom

\begin{abstract}
For real numbers $\alpha$ and $\beta$, let $\Bset(\alpha,\beta)= \{\lfloor\alpha m+\beta\rfloor:m\in\N\}$. Wooley asked whether the exponential sum over $\Bset(\alpha,\beta)\cap[1,N]$ always has $L^1$ norm $O_{\alpha,\beta}(\log N)$. We give a negative answer. We construct fixed $\alpha>1$ and $\beta$, with $1/\alpha$ of finite Diophantine type, such that the $L^1$ norm is $\gg(\log N)^2$ for infinitely many $N$. Wooley also proved that this $L^1$ norm is $O_{\alpha,\beta}((\log(2N))^2)$ when $1/\alpha$ has finite Diophantine type. Our example shows that his upper bound is sharp.
\end{abstract}

\maketitle

\section{Introduction}
Throughout the paper, for convenience, we write $\e(x)=\exp(2\pi i x)$, $\|x\|=\min_{m\in\mathbb Z}|x-m|$. For real numbers $\alpha$ and $\beta$, let $\Bset(\alpha,\beta) =\{\lfloor\alpha m+\beta\rfloor:m\in\N\}$ with $\N=\{1,2,3,\ldots\}$. 

Beatty sequences have a long history, originating with a question about how the positive integers can be partitioned into two complementary sequences. A classical partition theorem, known already to Rayleigh and later posed as a problem by Beatty \cite{Rayleigh,Beatty}, states that if $a,b>1$ are irrational and $1/a+1/b=1$, then the two sequences $\Bset(a,0)$ and $\Bset(b,0)$ partition the positive integers. Although their terms are defined via the floor function, the two sequences interlock perfectly, with no gaps or overlaps. This simple fact is one reason Beatty sequences have remained a natural subject in the study of integer sequences and their distribution. 

The partition theorem also gives a simple identity for exponential sums. For every $N$, we have
$\sum_{n\le N,\,n\in\Bset(a,0)}\e(n\theta)
+\sum_{n\le N,\,n\in\Bset(b,0)}\e(n\theta)
= \sum_{n=1}^N\e(n\theta).$
The sum on the right has $L^1$ norm $\asymp\log N$ for $N\ge2$. However, this does not imply the same bound for either Beatty sum separately.

Wooley asked whether the same logarithmic upper bound holds for every fixed Beatty sequence \cite[p.~2760]{Oberwolfach}. More precisely, he asked whether, for every fixed pair $\alpha,\beta$,
\begin{equation}\label{eq:wooley-question}
\int_0^1
\left|
\sum_{\substack{1\le n\le N\\
n\in\Bset(\alpha,\beta)}}
\e(n\theta)
\right|\,d\theta
\ll_{\alpha,\beta}\log N
\qquad(N\ge2).
\end{equation}

Wooley proved \cite[Theorem~2.3]{Wooley} that, when $\alpha>0$ and $1/\alpha$ has finite Diophantine type,
\begin{equation}\label{eq:wooley-upper}
\int_0^1
\left|
\sum_{\substack{1\le n\le N\\
n\in\Bset(\alpha,\beta)}}
\e(n\theta)
\right|\,d\theta
\ll_{\alpha,\beta}(\log(2N))^2.
\end{equation}
We say that an irrational number $\gamma$ has finite Diophantine type if $\inf_{Q\in\N}Q^k\|Q\gamma\|>0$ for some positive integer $k$. Our main result is the following.

\Needspace{10\baselineskip}
\begin{theorem}\label{thm:main}
There are fixed real numbers $\alpha>1$ and $\beta=-\alpha/2$ such that $\gamma=1/\alpha$ satisfies
\begin{equation}\label{eq:type-two}
\inf_{Q\in\N}Q^2\|Q\gamma\|>0,
\end{equation}
and there is a sequence of positive integers $N_j\to\infty$ for which
\begin{equation}\label{eq:main-lower}
\int_0^1
\left|
\sum_{\substack{1\le n\le N_j\\
n\in\Bset(\alpha,\beta)}}
\e(n\theta)
\right|\,d\theta
\gg(\log N_j)^2.
\end{equation}
\end{theorem}

Thus, the proposed bound \eqref{eq:wooley-question} fails, and the upper bound \eqref{eq:wooley-upper} is sharp under the finite Diophantine type assumption.

Our proof starts with a binary sequence obtained by coding a suitable rational rotation $p/q$. This sequence has period $q$. Using Fourier analysis, we show that its exponential sum over the first $q^2$ terms has $L^1$ norm at least a constant times $(\log q)^2$. We then use continued fractions to construct a fixed irrational parameter of finite Diophantine type. For infinitely many $q$, the first $q^2$ terms of the irrational sequence agree exactly with those of the corresponding rational sequence. This gives the same lower bound for the irrational Beatty sequence.

\section{Preliminaries}

We first record the relation between Beatty sequences and mechanical words \cite{MorseHedlund}, together with the Fourier identities used in the proof. For $0<\rho<1$ and $n\in\mathbb Z$, define
\begin{equation}\label{eq:mechanical-word}
  w_n(\rho)
  =\left\lfloor(n+1)\rho+\frac12\right\rfloor
   -\left\lfloor n\rho+\frac12\right\rfloor.
\end{equation}
Since $0<\rho<1$, every $w_n(\rho)$ is either $0$ or $1$. The shift by $1/2$ will be important when we compare rational and irrational words.

For our choice of intercept, the mechanical word is exactly the indicator of the corresponding Beatty sequence. In particular, there is no approximation error in passing between the two descriptions.

\begin{lemma}\label{lem:beatty-mechanical}
Let $0<\gamma<1$ be irrational, and set $\alpha=1/\gamma$ and $\beta=-\alpha/2$. Then, for every $n\in\N$,
$$
  \ind_{\Bset(\alpha,\beta)}(n)=w_n(\gamma).
$$
\end{lemma}

\begin{proof}
The condition $n=\lfloor\alpha m+\beta\rfloor$ is equivalent to $n\gamma+1/2\le m<(n+1)\gamma+1/2$. This interval has length $\gamma<1$ and positive, nonintegral endpoints. Its number of integers is therefore the indicator on the left and also $\lfloor(n+1)\gamma+1/2\rfloor-\lfloor n\gamma+1/2\rfloor$, which equals $w_n(\gamma)$.
\end{proof}

For $M\ge1$, let $D_M(x)=\sum_{r=0}^{M-1}\e(rx)$ be the Dirichlet kernel. The geometric-series identity gives, for $x\notin\mathbb Z$,
\begin{equation}\label{eq:dirichlet-formula}
  |D_M(x)|=\frac{|\sin(\pi Mx)|}{|\sin(\pi x)|}.
\end{equation}
We also use the discrete Parseval identity. For any complex numbers
$c_1,\ldots,c_q$,
$$
  \frac1q\sum_{\ell=0}^{q-1}
  \left|\sum_{n=1}^q c_n\e\left(\frac{n\ell}{q}\right)\right|^2
  =\sum_{n=1}^q|c_n|^2.
$$
Both are standard Fourier identities; the latter follows at once from orthogonality modulo $q$ \cite[Chapter~I]{Katznelson}.

\section{Rational mechanical words and the lower bound}

Let $p,q$ be coprime positive integers with $p<q$ and $q$ odd. The word $w_n(p/q)$ is $q$-periodic and has exactly $p$ ones in each period.
Write
$$
  P_{p,q}(x)=\sum_{n=1}^q w_n(p/q)\e(nx),\qquad
  A_{p,q}(u)=\frac1q\sum_{\ell=0}^{q-1}
  \left|P_{p,q}\left(\frac{\ell+u}{q}\right)\right|
  \quad(0\le u\le1).
$$
Here $P_{p,q}$ is the Fourier polynomial of one period, and $A_{p,q}(u)$ averages its absolute value over a shifted grid. The next lemma shows that this average is logarithmically large at $u=0$ and changes little when $u$ is small.

\begin{lemma}\label{lem:block-samples}
Fix $\eta\in(0,1/2]$. There is $c_\eta>0$ such that, for all coprime positive integers $p,q$ with $q$ odd and $\eta\le p/q\le1-\eta$,
$$
  A_{p,q}(0)\ge c_\eta\log q.
$$
Moreover, for $0\le u\le1$,
$$
  A_{p,q}(u)\ge A_{p,q}(0)-2\pi u\sqrt q.
$$
\end{lemma}

\begin{proof}
Write $\rho=p/q$. By \eqref{eq:mechanical-word}, $w_n(\rho)=1$ precisely when $\{n\rho+1/2\}\in[1-\rho,1)$. Since $q$ is odd, the fractional parts $\{np/q+1/2\}$ run through the grid $\{(j+1/2)/q:0\le j<q\}$, whose points avoid both endpoints of this interval. Thus the support of $w_n(p/q)$, after multiplication by $p$ modulo $q$, consists of $p$ consecutive residues.

Let $\bar p$ be the inverse of $p$ modulo $q$. For a suitable integer $b$, the support can be written as
$n\equiv\bar p(b+t)\pmod q$ for $0\le t<p$. Hence
$$
  P_{p,q}\left(\frac\ell q\right)
  =\e\left(\frac{\ell\bar p b}{q}\right)
   D_p\left(\frac{\ell\bar p}{q}\right),
  \qquad
  A_{p,q}(0)=\frac1q\sum_{j=0}^{q-1}
   \left|D_p\left(\frac jq\right)\right|.
$$
This is the finite Fourier calculation behind the formula of Graham and O'Bryant \cite[Theorem~1.1]{GrahamOBryant}.

For $1\le j\le q/2$, \eqref{eq:dirichlet-formula} gives $q^{-1}|D_p(j/q)|\ge|\sin(\pi\rho j)|/(\pi j)$. The inequality
$|\sin a|+|\sin b|\ge|\sin(b-a)|$ then shows that each pair $j=2r-1,2r$ contributes at least 
$\sin(\pi\rho)/(2\pi r)\ge\sin(\pi\eta)/(2\pi r)$.
Summing these contributions gives
$$
  A_{p,q}(0)\ge\frac pq+
  \frac{\sin(\pi\eta)}{2\pi}
  \sum_{r=1}^{\lfloor q/4\rfloor}\frac1r
  \gg_\eta\log q.
$$
The term $p/q\ge\eta$ also covers the case $q=3$.

For the second assertion, the reverse triangle inequality, Cauchy--Schwarz, and discrete Parseval give
\begin{align*}
  A_{p,q}(0)-A_{p,q}(u)
  &\le\left(\frac1q\sum_{\ell=0}^{q-1}
    \left|P_{p,q}\left(\frac{\ell+u}{q}\right)
         -P_{p,q}\left(\frac\ell q\right)\right|^2\right)^{1/2}\\
  &=\left(\sum_{n=1}^q w_n(p/q)
       |\e(nu/q)-1|^2\right)^{1/2}
  \le 2\pi u\sqrt q,
\end{align*}
where $|\e(t)-1|\le2\pi|t|$ and $n\le q$.
\end{proof}

We shall need a short-range form of the familiar logarithmic lower bound for the Dirichlet kernel. We give the elementary argument because this estimate supplies the second logarithm in the final bound.

\begin{lemma}\label{lem:short-dirichlet}
There is an absolute constant $c>0$ such that, if $M\ge2$, $0<U\le1/2$, and $MU\ge2$, then
$$
  \int_0^U|D_M(u)|\,du\ge c\log(MU).
$$
\end{lemma}

\begin{proof}
Put $T=MU$. Since $|\sin(\pi u)|\le\pi u$, the change of variables $t=Mu$ gives
$$
  \int_0^U|D_M(u)|\,du
  \ge\frac1\pi\int_0^T\frac{|\sin(\pi t)|}{t}\,dt
  \ge\frac2{\pi^2}\sum_{r=1}^{\lfloor T\rfloor}\frac1r
  \ge\frac2{\pi^2}\log T.
$$
Here we used $\int_{r-1}^r|\sin(\pi t)|\,dt=2/\pi$ and $1/t\ge1/r$ on $[r-1,r]$.
\end{proof}

We now repeat one period of the rational mechanical word $q$ times. The resulting exponential sum factors into $P_{p,q}$ and a Dirichlet kernel. The two logarithmic estimates above then combine to give the lower bound that will be transferred to an irrational Beatty sequence.

\begin{proposition}\label{prop:rational-repeat}
Fix $\eta\in(0,1/2]$. Let $p,q$ be coprime positive integers with $q$ odd and $\eta\le p/q\le1-\eta$. For all sufficiently large $q$
in terms of $\eta$,
$$
  \int_0^1\left|\sum_{n=1}^{q^2}w_n(p/q)\e(n\theta)\right|\,d\theta
  \gg_\eta(\log q)^2.
$$
\end{proposition}

\begin{proof}
Write $n=s+rq$ with $1\le s\le q$ and $0\le r<q$; this includes $n=q^2$, for which $s=q$ and $r=q-1$. Periodicity gives
\begin{equation}\label{eq:factorization}
  \sum_{n=1}^{q^2}w_n(p/q)\e(n\theta)
  =P_{p,q}(\theta)D_q(q\theta).
\end{equation}
Partition $[0,1]$ into intervals of length $1/q$ and substitute $\theta=(\ell+u)/q$. Using the periodicity of $D_q$, we obtain
\begin{equation}\label{eq:integral-factorization}
  \int_0^1|P_{p,q}(\theta)D_q(q\theta)|\,d\theta
  =\int_0^1|D_q(u)|A_{p,q}(u)\,du.
\end{equation}

Set $U=q^{-1/2}$. For sufficiently large $q$, Lemma~\ref{lem:block-samples} gives $A_{p,q}(u)\ge c_\eta\log q-2\pi\ge(c_\eta/2)\log q$ whenever $0\le u\le U$. By Lemma~\ref{lem:short-dirichlet} and \eqref{eq:integral-factorization},
$$
  \int_0^1|P_{p,q}(\theta)D_q(q\theta)|\,d\theta
  \ge\frac{c_\eta}{2}\log q\int_0^U|D_q(u)|\,du
  \gg_\eta(\log q)^2.
$$
Together with \eqref{eq:factorization}, this proves the proposition.
\end{proof}

\section{Construction of the irrational parameter}

We now choose the irrational parameter $\gamma$. To transfer the rational lower bound, we need good rational approximations to $\gamma$. At the same time, $\gamma$ must have finite Diophantine type. We arrange both properties by controlling the partial quotients of its continued fraction. The same choice will ensure that the denominators we use are odd.

Let $\gamma=[0;a_1,a_2,\ldots]$, and write $p_n/q_n$ for its convergents, with $q_{-1}=0$ and $q_0=1$. Choose the partial quotients
recursively by
\begin{equation}\label{eq:partial-quotients}
  a_{2r+1}=8q_{2r},\qquad a_{2r+2}=2\qquad(r\ge0),
\end{equation}
where $q_{n+1}=a_{n+1}q_n+q_{n-1}$. We use the standard best-approximation and error estimates for continued fractions \cite[Chapter~X]{HardyWright}.

The large partial quotients at odd indices give particularly accurate convergents at even indices. The denominators grow slowly enough that $\gamma$ has finite Diophantine type. The parity of the denominators will be needed in the next lemma.

\begin{lemma}\label{lem:gamma-properties}
The number $\gamma$ is irrational and lies in $(0,1)$. For every $r\ge0$, the denominator $q_{2r}$ is odd and $q_{2r+1}$ is even. Moreover, $\inf_{Q\in\N}Q^2\|Q\gamma\|>0$. For $r\ge1$ and $p/q=p_{2r}/q_{2r}$, we have
\begin{equation}\label{eq:selected-approximation}
  \left|\gamma-\frac pq\right|<\frac1{8q^3}.
\end{equation}
\end{lemma}

\begin{proof}
The continued fraction is infinite, so $\gamma$ is irrational; since $a_1=8$, it satisfies $1/9<\gamma<1/8$. Every partial quotient is even. The recurrence for $q_n$ therefore gives $q_{n+1}\equiv q_{n-1}\pmod 2$, proving the parity assertions.

We have $a_{n+1}\le8q_n$, and hence $q_{n+1}\le9q_n^2$ for every $n\ge0$. If $Q\ge q_1$, choose $n\ge1$ with $q_n\le Q<q_{n+1}$. The best-approximation property of convergents gives
$$
  \|Q\gamma\|\ge\|q_n\gamma\|
  >\frac1{q_n+q_{n+1}}
  \ge\frac1{10q_n^2}\ge\frac1{10Q^2}.
$$
For $1\le Q<q_1$, the quantity $\|Q\gamma\|$ is positive because $\gamma$ is irrational. Thus $\inf_{Q\in\N}Q^2\|Q\gamma\|>0$.

Finally, if $n=2r$, then $a_{n+1}=8q_n$. The usual convergent estimate gives $|\gamma-p_n/q_n|<1/(a_{n+1}q_n^2)=1/(8q_n^3)$, as required.
\end{proof}

These selected convergents do more than approximate $\gamma$ closely. With our half-integer shift, the corresponding rational and irrational mechanical words agree term by term through the first $q^2$ positions. This exact agreement lets us use the rational lower bound without an additional error estimate.

\begin{lemma}\label{lem:exact-shadow}
Let $p/q=p_{2r}/q_{2r}$ with $r\ge1$. Then
\begin{equation}\label{eq:exact-shadow}
  w_n(\gamma)=w_n(p/q)\qquad(1\le n\le q^2).
\end{equation}
\end{lemma}

\begin{proof}
Because $q$ is odd, $2mp+q$ is odd for every integer $m$. Thus
$$
  \left\|\frac{mp}{q}+\frac12\right\|
  \ge\frac1{2q}.
$$
On the other hand, \eqref{eq:selected-approximation} shows that, for
$1\le m\le q^2+1$,
$$
  \left|m\gamma-\frac{mp}{q}\right|
  <\frac{q^2+1}{8q^3}<\frac1{4q}.
$$
This change is smaller than the distance of the rational argument to any integer. Therefore
$$
  \left\lfloor m\gamma+\frac12\right\rfloor
  =\left\lfloor\frac{mp}{q}+\frac12\right\rfloor
  \qquad(1\le m\le q^2+1).
$$
Subtracting the equalities at $m=n+1$ and $m=n$ proves \eqref{eq:exact-shadow}.
\end{proof}

\Needspace{11\baselineskip}
\section{Proof of the main theorem}

We now combine Proposition~\ref{prop:rational-repeat} with Lemma~\ref{lem:exact-shadow} to prove the main theorem.

\begin{proof}
Let $\gamma$ be defined by \eqref{eq:partial-quotients}, and set $\alpha=1/\gamma$ and $\beta=-\alpha/2$.
Lemma~\ref{lem:gamma-properties} gives $\alpha>1$ and \eqref{eq:type-two}.

For $r\ge1$, write $p/q=p_{2r}/q_{2r}$ and $N_r=q^2$. Since $p_{2r}/q_{2r}\to\gamma\in(0,1)$, there is a fixed
$\eta\in(0,1/2)$ such that $\eta\le p/q\le1-\eta$ for all sufficiently large $r$. By Lemmas~\ref{lem:beatty-mechanical} and~\ref{lem:exact-shadow},
$$
  \sum_{\substack{1\le n\le N_r\\ n\in\Bset(\alpha,\beta)}}
  \e(n\theta)
  =\sum_{n=1}^{q^2}w_n(p/q)\e(n\theta).
$$
Since $q$ is odd, Proposition~\ref{prop:rational-repeat} gives
$$
  \int_0^1
  \left|
    \sum_{\substack{1\le n\le N_r\\ n\in\Bset(\alpha,\beta)}}
    \e(n\theta)
  \right|\,d\theta
  \gg_\eta(\log q)^2
  =\frac14(\log N_r)^2.
$$
As $q_{2r}\to\infty$, we have $N_r\to\infty$. Discarding finitely many initial terms gives \eqref{eq:main-lower}.
\end{proof}

\begin{remark}\label{rem:intercept}
The choice $\beta=-\alpha/2$ gives the half-integer shift in \eqref{eq:mechanical-word}. When $q$ is odd, the numbers $mp/q+1/2$ stay at distance at least $1/(2q)$ from the integers. This separation allows the exact agreement in Lemma~\ref{lem:exact-shadow}.
\end{remark}

\begin{remark}\label{rem:p-greater-one}
By Wooley's theorem \cite[Theorem~2.3]{Wooley}, for every fixed $s>1$, the same Beatty sequence satisfies
$$
  \int_0^1
  \left|
    \sum_{\substack{1\le n\le N\\ n\in\Bset(\alpha,\beta)}}
    \e(n\theta)
  \right|^s\,d\theta
  \ll_{\alpha,\beta,s}N^{s-1}.
$$
Thus our example rules out the proposed logarithmic upper bound at the $L^1$ endpoint, while the higher-moment estimates remain valid.
\end{remark}

\section*{Acknowledgments}

The author is very grateful to Niclas Technau for several helpful discussions about this problem. In particular, his questions about the factorization \eqref{eq:factorization} made the author realize the importance of checking the identities and indexing with greater care. His comments and suggestions were very useful during the development of the paper.

\enlargethispage{2\baselineskip}
\section*{Use of artificial intelligence}

ChatGPT was used during the research and preparation of this paper. In particular, ChatGPT 5.6 brought the work of Graham and O'Bryant \cite{GrahamOBryant} on the discrete Fourier transform of rational Beatty sets to the author's attention. Their Fourier formula provided the key starting point for our approach based on rational periodic sequences. Another decisive step was to construct a fixed irrational parameter of finite Diophantine type. For infinitely many $q$, its mechanical word agrees exactly with a rational mechanical word on the first $q^2$ terms. This construction was developed by the author through discussions with ChatGPT.

ChatGPT also assisted in checking calculations, identifying errors in earlier drafts, and improving the exposition and LaTeX source. The author is responsible for all mathematical arguments and results in the paper.


\begin{thebibliography}{9}

\bibitem{Beatty}
S. Beatty,
\emph{Problem 3173},
Amer. Math. Monthly \textbf{33} (1926), no.~3, 159,
\href{https://doi.org/10.2307/2300153}
{doi:10.2307/2300153}.

\bibitem{GrahamOBryant}
R. L. Graham and K. O'Bryant,
\emph{A discrete Fourier kernel and Fraenkel's tiling conjecture},
Acta Arith. \textbf{118} (2005), no.~3, 283--304,
\href{https://doi.org/10.4064/aa118-3-4}
{doi:10.4064/aa118-3-4}.

\bibitem{HardyWright}
G. H. Hardy and E. M. Wright,
\emph{An Introduction to the Theory of Numbers},
sixth edition, Oxford University Press, Oxford, 2008.

\bibitem{Katznelson}
Y. Katznelson,
\emph{An Introduction to Harmonic Analysis},
third edition, Cambridge University Press,
Cambridge, 2004.

\bibitem{MorseHedlund}
M. Morse and G. A. Hedlund,
\emph{Symbolic dynamics II. Sturmian trajectories},
Amer. J. Math. \textbf{62} (1940), no.~1, 1--42,
\href{https://doi.org/10.2307/2371431}
{doi:10.2307/2371431}.

\bibitem{Oberwolfach}
K. Matom\"aki, J. Maynard, K. Soundararajan,
and T. D. Wooley (organizers),
\emph{Analytic Number Theory},
Oberwolfach Rep. \textbf{22} (2025), no.~4,
2701--2764,
\href{https://doi.org/10.4171/OWR/2025/51}
{doi:10.4171/OWR/2025/51}.

\bibitem{Rayleigh}
J. W. Strutt (Lord Rayleigh),
\emph{The Theory of Sound},
vol.~I, second edition,
Macmillan, London, 1894, p.~123.

\bibitem{Wooley}
T. D. Wooley,
\emph{Subconvex $L^p$-sets, Weyl's inequality
and equidistribution},
J. Inst. Math. Jussieu \textbf{25} (2026),
no.~5, 2919--2959,
\href{https://doi.org/10.1017/S1474748026101935}
{doi:10.1017/S1474748026101935}.

\end{thebibliography}
\end{document}